\documentclass[a4paper, 12pt]{article}
\usepackage{amssymb, amsmath, amsfonts}
\usepackage{amsthm,amscd,amsfonts,latexsym,color,epsfig}
\begin{document}

\begin{center}
{\bf On a non-local problem for parabolic-hyperbolic equation with three lines of type changing}\\
Karimov E.T., Sotvoldiev A.I.\\
E-mail: erkinjon@gmail.com,akm0111@inbox.ru\\
\emph{Institute of Mathematics, National University of Uzbekistan named after Mirzo Ulughbek (Tashkent, Uzbekistan)}\\
\end{center}

\bigskip


\bigskip

\textbf{Abstract. }In the present work we investigate a boundary problem with non-local conditions, connecting values of seeking function on various characteristics for parabolic-hyperbolic equation with three lines of type changing. The considered problem is equivalently reduced to the system of Volterra integral equations of the second kind.

\bigskip

Consider an equation
$$
0 = \left\{ \begin{gathered}
  {u_{xx}} - {u_y},\,\,\,\,\,\,\left( {x,y} \right) \in {\Omega _0}, \hfill \\
  {u_{xx}} - {u_{yy}},\,\,\,\,\,\left( {x,y} \right) \in {\Omega _i}\,\left( {i = \overline {1,3} } \right) \hfill \\
\end{gathered}  \right.\eqno (1)
$$
in the domain $\Omega  = {\Omega _0} \cup {\Omega _1} \cup {\Omega _2} \cup {\Omega _3} \cup AB \cup A{A_0} \cup B{B_0}$.
\begin{figure}[h]
\begin{center}
\includegraphics[width=0.4\textwidth]{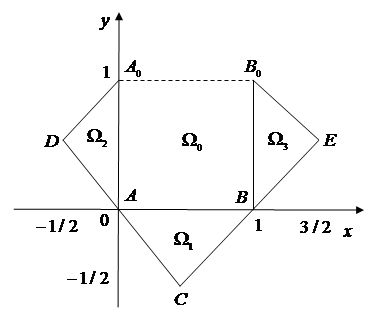}\\
\caption{Domain $\Omega$}
\end{center}
\end{figure}

\textbf{Problem AS.} Find a regular solution of the equation (1) in the domain $\Omega $, satisfying the following conditions:
$$
{a_1}\left( t \right)u\left( { - t,t} \right) + {a_2}\left( t \right)u\left( {t, - t} \right) = {a_3}\left( t \right),\,\,0 \leqslant t \leqslant \frac{1}{2},\eqno (2)
$$
$$
{b_1}\left( t \right)u\left( {t,t - 1} \right) + {b_2}\left( t \right)u\left( {2 - t,1 - t} \right) = {b_3}\left( t \right),\,\,\frac{1}{2} \leqslant t \leqslant 1,\eqno (3)
$$
$$
{c_1}\left( t \right)\left( {{u_x} + {u_y}} \right)\left( {t - 1,t} \right) + {c_2}\left( t \right)\left( {{u_x} - {u_y}} \right)\left( {2 - t,t} \right) = {c_3}\left( t \right),\,\,\frac{1}{2} < t < 1.\eqno (4)
$$
Here ${a_i}\left( t \right),{b_i}\left( t \right),{c_i}\left( t \right)\,\left( {i = \overline {1,3} } \right)$ are given functions, such that
$$
\begin{gathered}
  {a_1}\left( 0 \right) + {a_2}\left( 0 \right) \ne 0,\,{b_1}\left( 1 \right) + {b_2}\left( 1 \right) \ne 0,\,a_1^2\left( t \right) + a_2^2\left( t \right) > 0,\,b_1^2\left( t \right) + b_2^2\left( t \right) > 0, \hfill \\
  c_1^2\left( t \right) + c_2^2\left( t \right) > 0,\,a_1^2 + b_2^2 > 0,\,a_2^2 + b_1^2 > 0. \hfill \\
\end{gathered}
$$

Note, boundary problems for parabolic-hyperbolic equations with two lines of type changing were investigated in the works [1-4], and with three lines of type changing in the papers [5-6]. Distinctive side of the present work is non-local condition, which connect values of seeking function on various characteristics. It makes very difficult the reduction of the considered problem to the system of integral equations and we need special algorithm for solving this problem.

In the domain ${\Omega _1}$ solution of the Cauchy problem with initial data $u\left( {x,0} \right) = {\tau _1}\left( x \right)$, ${u_y}\left( {x,0} \right) = {\nu _1}\left( x \right)$ can be represented as
$$
2u\left( {x,y} \right) = {\tau _1}\left( {x + y} \right) + {\tau _1}\left( {x - y} \right) + \int\limits_{x - y}^{x + y} {{\nu _1}\left( z \right)dz} .\eqno (5)
$$
Assuming in condition (2)
$$
u\left( { - t,t} \right) = {\varphi _1}\left( t \right),\,\,\,0 \leqslant t \leqslant \frac{1}{2},\eqno (6)
$$
as given, from (5) we find
$$
{\tau '_1}\left( t \right) = {\nu _1}\left( t \right) + {\left( {\frac{{2\left[ {{a_3}\left( {\frac{t}{2}} \right) - {a_1}\left( {\frac{t}{2}} \right){\varphi _1}\left( {\frac{t}{2}} \right)} \right]}}{{{a_2}\left( {\frac{t}{2}} \right)}}} \right)^\prime },\,\,\,\,0 < t < 1.\eqno (7)
$$
In condition (3) introduce designation
$$
u\left( {2 - t,1 - t} \right) = {\varphi _2}\left( t \right),\,\,\,\frac{1}{2} \leqslant t \leqslant 1,\eqno (8)
$$
and from (5) we get
$$
{\tau '_1}\left( t \right) =  - {\nu _1}\left( t \right) + {\left( {\frac{{2\left[ {{b_3}\left( {\frac{{t + 1}}{2}} \right) - {b_2}\left( {\frac{{t + 1}}{2}} \right){\varphi _2}\left( {\frac{{t + 1}}{2}} \right)} \right]}}{{{b_1}\left( {\frac{{t + 1}}{2}} \right)}}} \right)^\prime },\,\,\,\,0 < t < 1.\eqno	(9)
$$
From (7) and (9) it follows that
$$
{\tau '_1}\left( t \right) = {\left( {\frac{{{a_3}\left( {\frac{t}{2}} \right) - {a_1}\left( {\frac{t}{2}} \right){\varphi _1}\left( {\frac{t}{2}} \right)}}{{{a_2}\left( {\frac{t}{2}} \right)}}} \right)^\prime } + {\left( {\frac{{{b_3}\left( {\frac{{t + 1}}{2}} \right) - {b_2}\left( {\frac{{t + 1}}{2}} \right){\varphi _2}\left( {\frac{{t + 1}}{2}} \right)}}{{{b_1}\left( {\frac{{t + 1}}{2}} \right)}}} \right)^\prime },\,\,\,\,0 < t < 1.\eqno (10)
$$

Solution of the Cauchy problem in the domain ${\Omega _2}$ with given data $u\left( {0,y} \right) = {\tau _2}\left( y \right)$, ${u_x}\left( {0,y} \right) = {\nu _2}\left( y \right)$ we write as follows
$$
2u\left( {x,y} \right) = {\tau _2}\left( {y + x} \right) + {\tau _2}\left( {y - x} \right) + \int\limits_{y - x}^{y + x} {{\nu _2}\left( z \right)dz} .\eqno	(11)
$$
Considering (6) from (11) we obtain
$$
{\tau '_2}\left( t \right) = {\nu _2}\left( t \right) + {\varphi '_1}\left( {\frac{t}{2}} \right),\,\,0 < t < 1.\eqno (12)
$$
In condition (4) introduce another designation
$$
\left( {{u_x} - {u_y}} \right)\left( {2 - t,t} \right) = {\varphi _3}\left( t \right),\,\,\frac{1}{2} < t < 1.\eqno (13)
$$
Then from (11) we get
$$
\frac{{{c_3}\left( {\frac{{t + 1}}{2}} \right) - {c_2}\left( {\frac{{t + 1}}{2}} \right){\varphi _3}\left( {\frac{{t + 1}}{2}} \right)}}{{{c_1}\left( {\frac{{t + 1}}{2}} \right)}} = {\tau '_2}\left( t \right) + {\nu _2}\left( t \right),\,\,0 < t < 1.\eqno (14)
$$
From (12) and (14) we deduce
$$
2{\tau '_2}\left( t \right) = {\varphi '_1}\left( {\frac{t}{2}} \right) + \frac{{{c_3}\left( {\frac{{t + 1}}{2}} \right) - {c_2}\left( {\frac{{t + 1}}{2}} \right){\varphi _3}\left( {\frac{{t + 1}}{2}} \right)}}{{{c_1}\left( {\frac{{t + 1}}{2}} \right)}},\,\,0 < t < 1.\eqno (15)
$$

Solution of the Cauchy problem with data $u\left( {1,y} \right) = {\tau _3}\left( y \right),\,{u_x}\left( {1,y} \right) = {\nu _3}\left( y \right)$ â in the domain ${\Omega _3}$ has a form
$$
2u\left( {x,y} \right) = {\tau _3}\left( {y + x - 1} \right) + {\tau _2}\left( {y - x + 1} \right) + \int\limits_{y - x + 1}^{y + x - 1} {{\nu _3}\left( z \right)dz}.\eqno (16)
$$
Using (8) and (13) from (16), after some evaluations one can get
$$
2{\tau '_3}\left( t \right) =  - {\varphi '_2}\left( {\frac{{2 - t}}{2}} \right) - {\varphi _3}\left( {\frac{{t + 1}}{2}} \right),\,\,0 < t < 1.\eqno (17)
$$
Further, from the equation (1) we pass to the limit at $y \to  + 0$ and considering (7) we find
$$
{\tau ''_1}\left( t \right) - {\tau '_1}\left( t \right) =  - {\left( {\frac{{2\left[ {{a_3}\left( {\frac{t}{2}} \right) - {a_1}\left( {\frac{t}{2}} \right){\varphi _1}\left( {\frac{t}{2}} \right)} \right]}}{{{a_2}\left( {\frac{t}{2}} \right)}}} \right)^\prime }.\eqno (18)
$$
Solution of the equation (18) together with conditions
$$
{\tau _1}\left( 0 \right) = \frac{{{a_3}\left( 0 \right)}}{{{a_1}\left( 0 \right) + {a_2}\left( 0 \right)}},\,\,\,{\tau _1}\left( 1 \right) = \frac{{{b_3}\left( 1 \right)}}{{{b_1}\left( 1 \right) + {b_2}\left( 1 \right)}},\eqno (19)
$$
which reduced from (2) and (3), can be represented as
$$
\begin{gathered}
  {\tau _1}\left( x \right) = \frac{{{a_3}\left( 0 \right)}}{{{a_1}\left( 0 \right) + {a_2}\left( 0 \right)}} + x\left[ {\frac{{{b_3}\left( 1 \right)}}{{{b_1}\left( 1 \right) + {b_2}\left( 1 \right)}} - \frac{{{a_3}\left( 0 \right)}}{{{a_1}\left( 0 \right) + {a_2}\left( 0 \right)}}} \right] +  \hfill \\
   + \int\limits_0^1 {G\left( {x,t} \right)} \left[ {\frac{{{b_3}\left( 1 \right)}}{{{b_1}\left( 1 \right) + {b_2}\left( 1 \right)}} - \frac{{{a_3}\left( 0 \right)}}{{{a_1}\left( 0 \right) + {a_2}\left( 0 \right)}}} \right]dt -  \hfill \\
   - \int\limits_0^1 {G\left( {x,t} \right)} {\left( {\frac{{2\left[ {{a_3}\left( {\frac{t}{2}} \right) - {a_1}\left( {\frac{t}{2}} \right){\varphi _1}\left( {\frac{t}{2}} \right)} \right]}}{{{a_2}\left( {\frac{t}{2}} \right)}}} \right)^\prime }dt,\,\,\,\,0 \leqslant x \leqslant 1, \hfill \\
\end{gathered}
\eqno 	(20)
$$
where $G\left( {x,t} \right)$ is Green's function of the problem (18)-(19).

Continuing to assume the function ${\varphi _1}$ as known, using the formula (10) we represent function ${\varphi _2}$ via ${\varphi _1}$. Then using the solutionf of the first boundary problem for the equation (1) in the domain ${\Omega _0}$ and functional relations between functions ${\tau _j}$ and ${\nu _j}$ $\left( {j = 2,3} \right)$, we get the following:
$$
\begin{gathered}
  {{\tau '}_2}\left( y \right) = \int\limits_0^y {{{\tau '}_3}\left( \eta  \right)N\left( {0,y,1,\eta } \right)d\eta }  - \int\limits_0^y {{{\tau '}_2}\left( \eta  \right)N\left( {0,y,0,\eta } \right)d\eta }  + {F_1}\left( y \right), \hfill \\
  {{\tau '}_3}\left( y \right) = \int\limits_0^y {{{\tau '}_3}\left( \eta  \right)N\left( {1,y,1,\eta } \right)d\eta }  - \int\limits_0^y {{{\tau '}_2}\left( \eta  \right)N\left( {1,y,0,\eta } \right)d\eta }  + {F_2}\left( y \right), \hfill \\
\end{gathered} \eqno (21)
$$
where
$$
\begin{gathered}
  {F_1}\left( y \right) = \int\limits_0^1 {{\tau _1}\left( \xi  \right){{\overline G }_x}\left( {o,y,\xi ,0} \right)d\xi }  - \frac{{{a_3}\left( 0 \right)}}{{{a_1}\left( 0 \right) + {a_2}\left( 0 \right)}}N\left( {0,y,0,0} \right) +  \hfill \\
   + \frac{{{b_3}\left( 1 \right)}}{{{b_1}\left( 1 \right) + {b_2}\left( 1 \right)}}N\left( {0,y,1,0} \right) + {{\varphi '}_1}\left( {\frac{y}{2}} \right), \hfill \\
\end{gathered}
$$
$$
\begin{gathered}
  {F_2}\left( y \right) = \int\limits_0^1 {{\tau _1}\left( \xi  \right){{\overline G }_x}\left( {1,y,\xi ,0} \right)d\xi }  - \frac{{{a_3}\left( 0 \right)}}{{{a_1}\left( 0 \right) + {a_2}\left( 0 \right)}}N\left( {1,y,0,0} \right) +  \hfill \\
   + \frac{{{b_3}\left( 1 \right)}}{{{b_1}\left( 1 \right) + {b_2}\left( 1 \right)}}N\left( {1,y,1,0} \right) - {\varphi _3}\left( {\frac{{y + 1}}{2}} \right), \hfill \\
\end{gathered}
$$
$$
\overline G \left( {x,y,\xi ,\eta } \right) = \frac{1}{{2\sqrt {\pi \left( {y - \eta } \right)} }}\sum\limits_{n =  - \infty }^\infty  {\left[ {{e^{ - \frac{{{{\left( {x - \xi  + 2n} \right)}^2}}}{{4\left( {y - \eta } \right)}}}} - {e^{ - \frac{{{{\left( {x + \xi  + 2n} \right)}^2}}}{{4\left( {y - \eta } \right)}}}}} \right]}
$$
is Green's function of the first boundary problem,
$$
N\left( {x,y,\xi ,\eta } \right) = \frac{1}{{2\sqrt {\pi \left( {y - \eta } \right)} }}\sum\limits_{n =  - \infty }^\infty  {\left[ {{e^{ - \frac{{{{\left( {x - \xi  + 2n} \right)}^2}}}{{4\left( {y - \eta } \right)}}}} + {e^{ - \frac{{{{\left( {x + \xi  + 2n} \right)}^2}}}{{4\left( {y - \eta } \right)}}}}} \right]}.
$$

From the first equation of  (21) we represent function ${\varphi _3}$ via ${\varphi _1}$ and further, from the second equation of (21) we find the function ${\varphi _1}$.

After the finding function ${\varphi _1}$, using appropriate formulas we find functions ${\varphi _2}$,${\varphi _3}$, ${\tau _i}$, ${\nu _i}$, $\left( {i = \overline {1,3} } \right)$. Solution of the problem AS can be established in the domain  ${\Omega _0}$ as a solution of the first boundary problem, and in the domains ${\Omega _i}\,\left( {i = \overline {1,3} } \right)\,$ as a solution of the Cauchy problem.

\textbf{Theorem.} If functions $a_i,\,b_i,\,c_i$ are continuously differentiable on the segment, and have continuous second order derivatives on interval, where they given, then the problem AS have the unique regular solution.

\smallskip

\centerline{\bf References}
\begin{enumerate}
\item {\it Egamberdiev U.} Boundary problems for mixed parabolic-hyperbolic equation with two lines of type changing. PhD thesis, Tashkent, 1984.\\
\item{\it Abdullaev A.S.} On some boundary problems for mixed parabolic-hyperbolic type equations// Equations of mixed type and problem with free boundary. Tashkent: Fan, 1987, pp. 71-82.\\
\item {\it Eleev V.A., Lesev V.N.} On two boundary problems for mixed type equations with perpendicular lines of type changing// Vladikavkaz math.journ. 2001. Vol. 3. Vyp.  4, pp.9-22.\\
\item{\it Nakusheva V.A.} First boundary problem for mixed type equation in a characteristic polygon// Dokl.AMAN, 2012. Vol.14, No 1, pp.58-65.\\
\item{\it Berdyshev A.S., Rakhmatullaeva N.A.} Nonlocal problems with special gluing for a parabolic-hyperbolic equation. "Further Progress in Analysis". Proceedings of the 6th ISAAC Congress. Ankara, Turkey, 13-18 August, 2007, pp. 727-734.\\
\item{\it Berdyshev A.S., Rakhmatullaeva N.A.} Non-local problems for parabolic-hyperbolic equations with deviation from the characteristics and three type-changing lines //Electronic Journal of Differential Equations. Vol. (2011) 2011, No 7, pp.1-6.\\
\end{enumerate}
\end{document}